\documentclass[final]{elsarticle}
\makeatletter
\def\ps@pprintTitle{%
 \let\@oddhead\@empty
 \let\@evenhead\@empty
 \def\@oddfoot{\centerline{\thepage}}%
 \let\@evenfoot\@oddfoot}
\makeatother

\usepackage{techRepsPDF}
\usepackage{myalg}
\usepackage{url}

\usepackage{amsthm}
\newtheorem{property}{Property}
\usepackage{amssymb}
\usepackage[nodots]{numcompress}

\usepackage[cmex10]{amsmath}
\usepackage{graphicx}
\usepackage{lscape}

\usepackage{enumitem}
\newlist{cenumerate}{enumerate}{1}
\setlist[cenumerate]{label=C-\arabic*}

\newcommand{\minim}{\mbox{minimize }} 
\newcommand{\st}{\mbox{subject to: }} 
\newcommand{\son}{\mathsf{on}}
\newcommand{\soff}{\mathsf{off}}
\newcommand{\xon}{x^{\son}} 
\newcommand{\xoff}{x^{\soff}} 
\newcommand{\add}{\mathsf{add}}

\newcommand{\hot}{\mathsf{hot}}
\newcommand{\cold}{\mathsf{cold}}
\newcommand{\prev}{\mathsf{prev}}
\newcommand{\Pmin}{P^{\mathsf{min}}} 
\newcommand{\Pmax}{P^{\mathsf{max}}} 
\newcommand{\for}{\mathop{\text{ for }}}
\newcommand{\ie}{\emph{i.e.}}
\newcommand{\eg}{\emph{e.g.}}

\trtitle{Unit commitment with valve-point loading effect}
\trnumber{DCC-2014-05}
\trauthor{
\normalsize 
Jo\~ao Pedro Pedroso\\
\textit{INESC Porto and Faculdade de Ci\^encias, Universidade do Porto, Portugal}\\
\texttt{jpp@fc.up.pt}\\
~\\
Mikio Kubo\\
\textit{Tokyo University of Marine Science and Technology, Japan}\\
\texttt{kubo@kaiyodai.ac.jp}\\
~\\
Ana Viana\\
\textit{INESC Porto and Instituto Superior de Engenharia, Instituto Polit\'ecnico do Porto, Portugal}\\
\texttt{aviana@inescporto.pt}
}

\begin{document}
\mkcoverpage
\begin{frontmatter}

\title{Unit commitment with valve-point loading effect}

\author{Jo\~ao Pedro Pedroso}
\ead{jpp@fc.up.pt}
\address{INESC Porto and Faculdade de Ci\^encias, Universidade do Porto, Portugal}

\author{Mikio Kubo}
\ead{kubo@kaiyodai.ac.jp}
\address{Tokyo University of Marine Science and Technology, Japan}

\author{Ana Viana}
\ead{aviana@inescporto.pt}
\address{INESC Porto and \\Instituto Superior de Engenharia, Instituto Polit\'ecnico do Porto, Portugal}

\date{April 2012}

\begin{abstract}
Valve-point loading affects the input-output characteristics of generating units, bringing the fuel costs nonlinear and nonsmooth.  This has been considered in the solution of load dispatch problems, but not in the planning phase of unit commitment.  

This paper presents a mathematical optimization model for the thermal unit commitment problem considering valve-point loading.  The formulation is based on a careful linearization of the fuel cost function, which is modeled with great detail on power regions being used in the current solution, and roughly on other regions.

A set of benchmark instances for this problem is used for analyzing the method, with recourse to a general-purpose mixed-integer optimization solver.
\end{abstract}

\begin{keyword}
Unit Commitment\sep
Load Dispatch\sep
Combinatorial Optimization\sep
Mixed-integer Programming
\end{keyword}
\end{frontmatter}

\section{Introduction}

The unit commitment problem (UCP) consists of deciding which power generating units must be committed/decommitted over a planning horizon, usually lasting from 1 day to 2 weeks and split into periods of one hour.  The production levels at which units must operate (pre-dispatch) must also be determined to optimize a cost function which includes both fixed and variable costs.  The committed units must satisfy the forecasted system load and reserve requirements, as well as a large set of technological constraints.  This problem has great practical relevance, due to the savings that can be achieved with an optimized schedule.

Simple approaches for this problem model the cost as a linear function, which provides a rather rough approximation; accuracy can be improved by switching to a quadratic function.  For these models general-purpose mixed-integer convex optimization solvers can be used to find optimum solutions for realistic, large instances.  Nonetheless, much better approximations of the true costs can be obtained if the valve-point loading effect is taken into account;  however models become much more difficult to optimize.  

Valve-point loading affects the input-output characteristics of generating units, making the fuel costs nonlinear and nonsmooth.  This has been considered in the heuristic solution of load dispatch problems --- see, \eg, \cite{AlSumait2007720, Hemamalini2011868, Roy20134244, SrinivasaReddy2013342}, but there are no exact methods available in the literature.

Concerning the planning phase of unit commitment, even for convex models most of the research on the solution of the UCP focuses in heuristic methods.   Recently, however, improvements in the capabilities of mixed-integer programming solvers has encouraged the thorough exploitation of their capabilities~\cite{viana2013}.  Extensive surveys of different optimization techniques and modeling issues are provided in \cite{sen1998, PAD04, YAM04}, but the valve-point loading effect is not considered in any of the exact methods proposed.  A mixed integer quadratically constrained model to solve unit commitment on real-life, large-scale power systems is presented in~\cite{Alvarez2012}; it details some features of the units, but not the valve-point loading effect.  Another approach, where a mixed integer linear formulation is used for modeling nonlinear output of generators and applied to unit commitment, is presented in~\cite{Simoglou2010}.  A process with some similarities to the method we propose, also involving the assessment of lower and upper bounds to a nonlinear function and iteratively converging to the optimum, has been presented in~\cite{Lima2013} and applied to a related problem: that of optimizing short-term hydro scheduling.

The main contribution of this paper is a method based on a careful linearization of the fuel cost function, which is modeled with great detail on power regions being used in the current solution and roughly on other regions; this formulation is used in an iterative process that converges to the optimum.  The method may be used both for load dispatch (by setting the number of planning periods to one) and for unit commitment.  When the solution process stops, \eg~due to a limit on CPU usage, the method provides both an upper and a lower bound to the optimum; this information is usually very important in practice.


\section{Mathematical model}
\label{sec:model}

A concise, complete description of the unit commitment problem as a mathematical optimization model has been provided in~\cite{viana2013}.  When valve-point loading is taken into account that model can be used as is, except for the shape of the objective function: instead of a quadratic function, as considered in~\cite{kazarlis1996}, it is now a nonconvex, nonsmooth function, as has been proposed for load dispatching in~\cite{wood1984} (cited by~\cite{coelho2006}).

\begin{figure}[!htbp]
\centering
\resizebox{.75\linewidth}{!}{\input{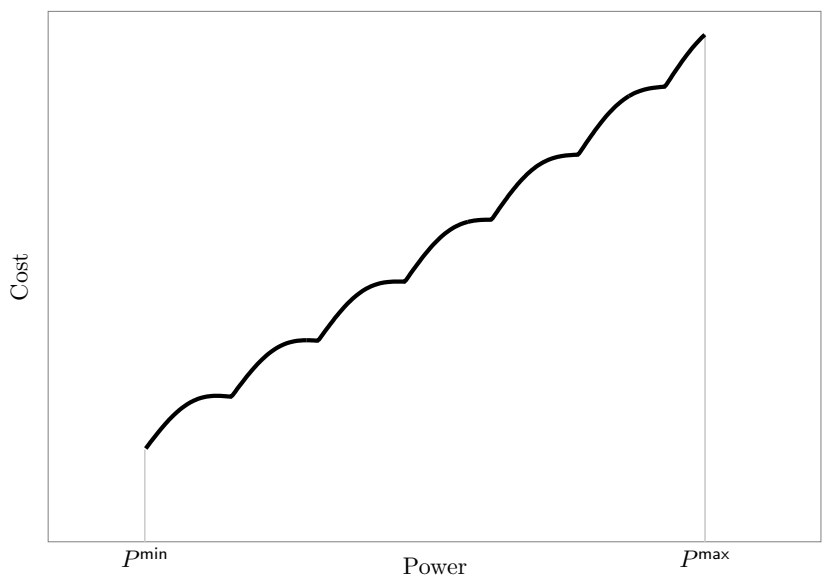}}
\caption{Typical shape of the fuel cost considering the valve-point loading effect.}
\label{fig:cost_shape}
\end{figure}

\subsection{Objective}
\label{sec:objective}

A commonly used function of the fuel costs in a generating unit in terms of the power produced $p$, taking into account the valve-point loading effect, is the following~\cite{wood1984}:
\begin{alignat}{27}
F(p) = a + b p + c p^2 + \left|e \sin\left(f (\Pmin-p)\right)\right|,  \label{eq:valve}
\end{alignat}
where $a,b,c,e,f$ are function's parameters (see Figure~\ref{fig:cost_shape}), as are the minimum and maximum operating powers $\Pmin$ and $\Pmax$.
The first three terms describe a quadratic function, as usually taken into account when valve-point loading is not considered; when it is, there is a periodic factor, as described in the rightmost term.

This function being nonlinear, nonconvex, and nonsmooth, methods for optimizing it under constraints on load limit and demand satisfaction --- \ie, in an economic load dispatch setting --- usually involve specifically designed heuristic approaches, many of them based on evolutionary algorithms~\cite{SrinivasaReddy2013342, Coelho20102580}.

\begin{figure}[!htbp]
\centering
\resizebox{.49\linewidth}{!}{\input{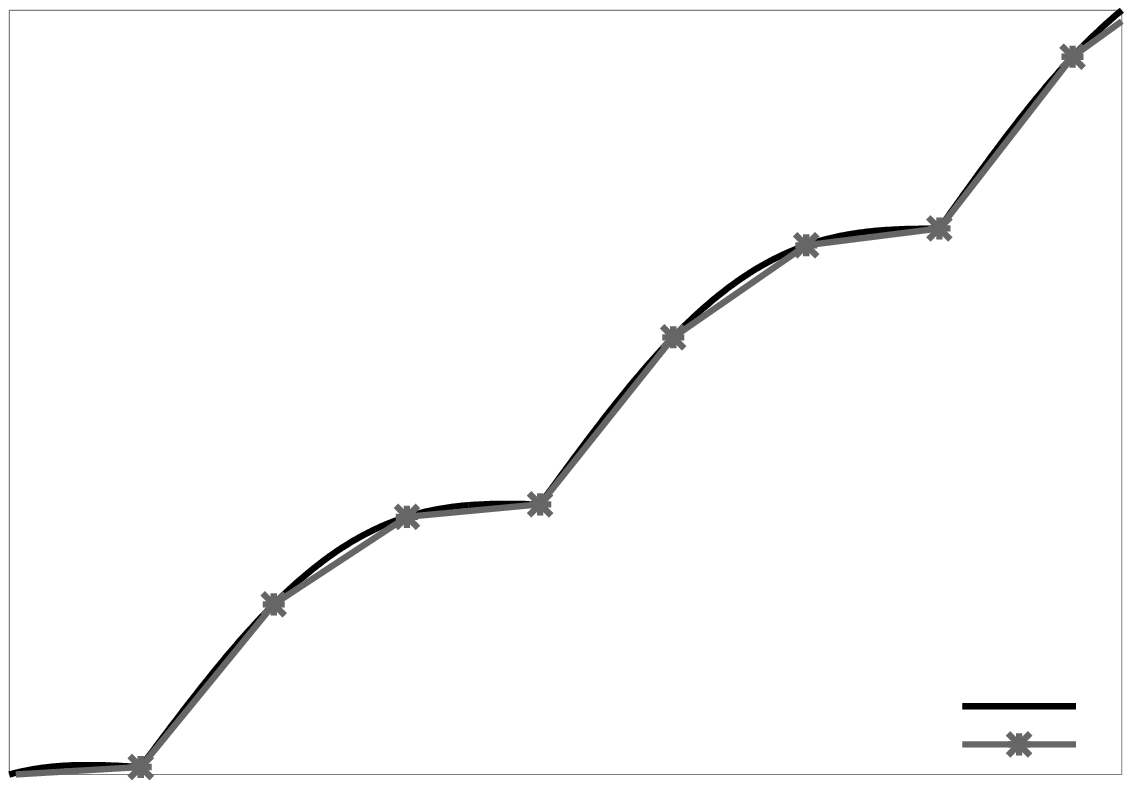}}
\resizebox{.49\linewidth}{!}{\input{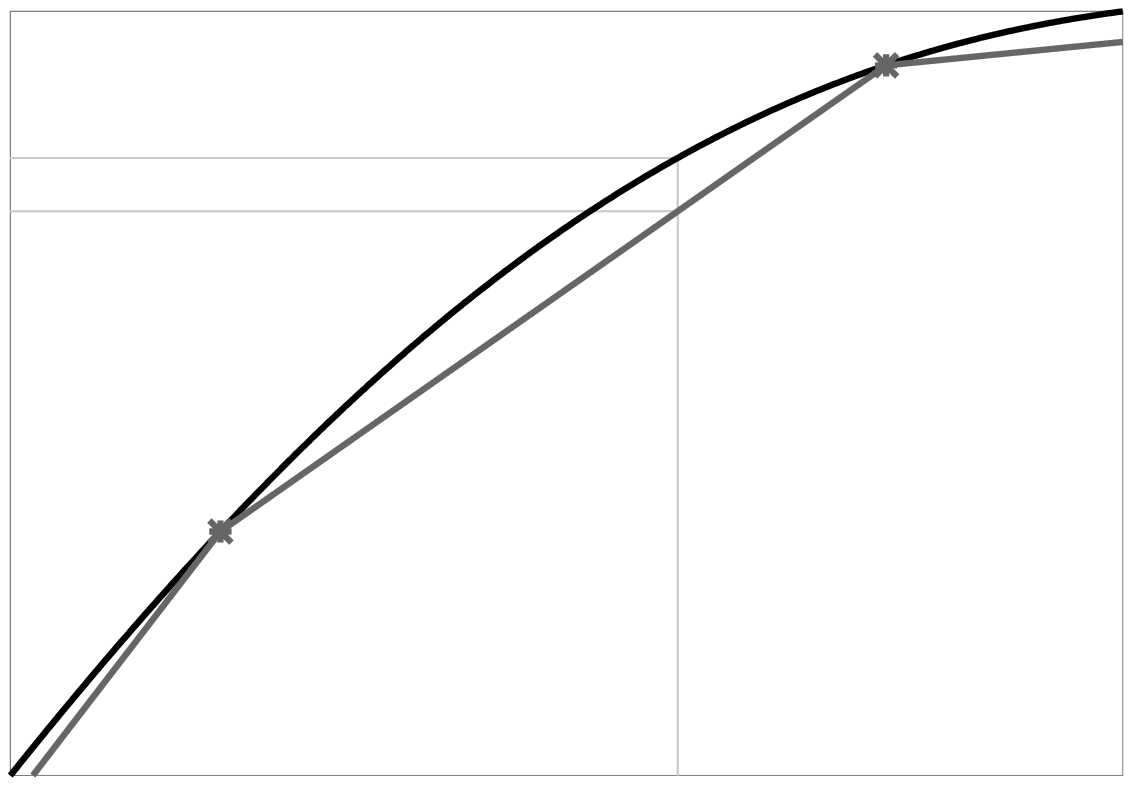}}
\caption{Fuel cost and its representation with a piecewise-linear function (left).  Upper bound (UB, exact function evaluation) and lower bound (LB, piecewise-linear approximation) at a given production level $p$ (right).}
\label{fig:cost_lowerA}
\label{fig:cost_bounds}
\end{figure}

There is, however, a way of tackling this problem exactly (\ie, with arbitrarily low error).  Let us first notice that the valve points --- where there is a discontinuity --- limit areas where the function is concave.  Thus, if we evaluate the function into several points and use them as breakpoints for joining the two valve points, we obtain a linear function that is never larger than the true function (see Figure~\ref{fig:cost_lowerA}).  Given a power production level this function provides, thus, a lower bound to the objective value; on the other hand, the evaluation of the true function provides a trivial upper bound.  For any solution satisfying all the problem's constraints, a lower and an upper bound to possible values of the objective can therefore be obtained by the piecewise-linear lower approximation and the true function, respectively.

\begin{figure}[!htbp]
\centering
\resizebox{.49\linewidth}{!}{\input{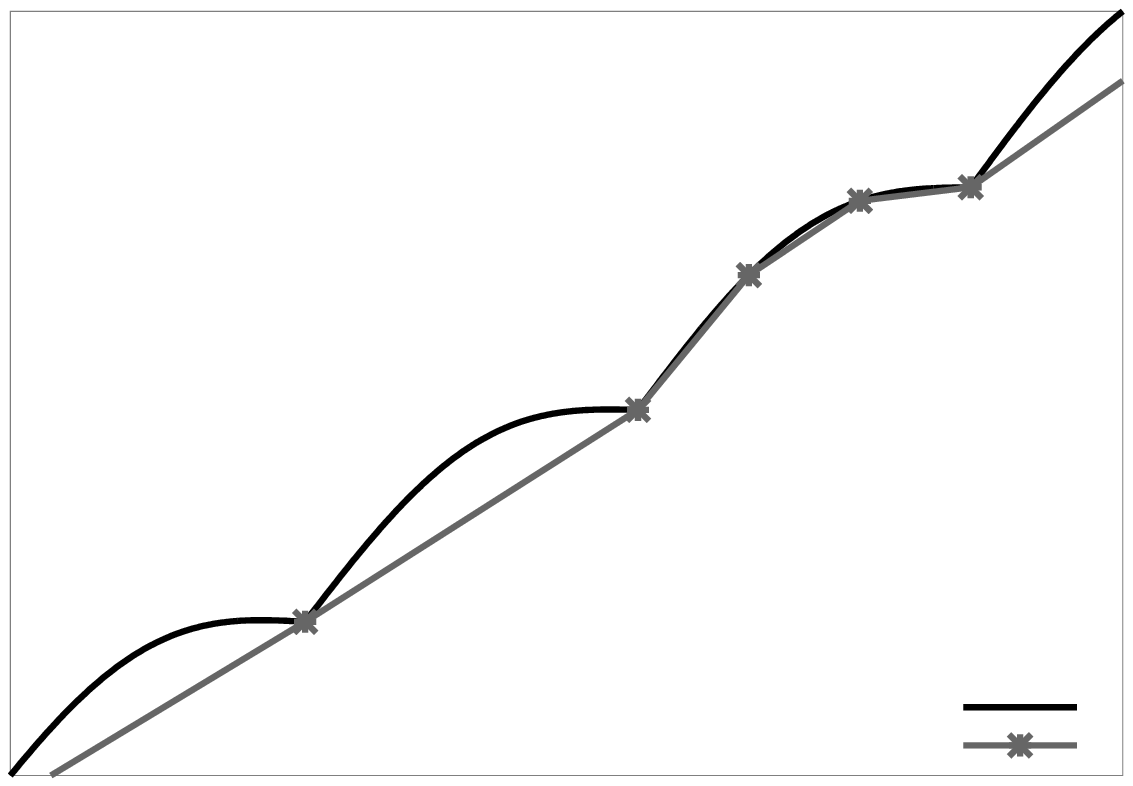}}
\resizebox{.49\linewidth}{!}{\input{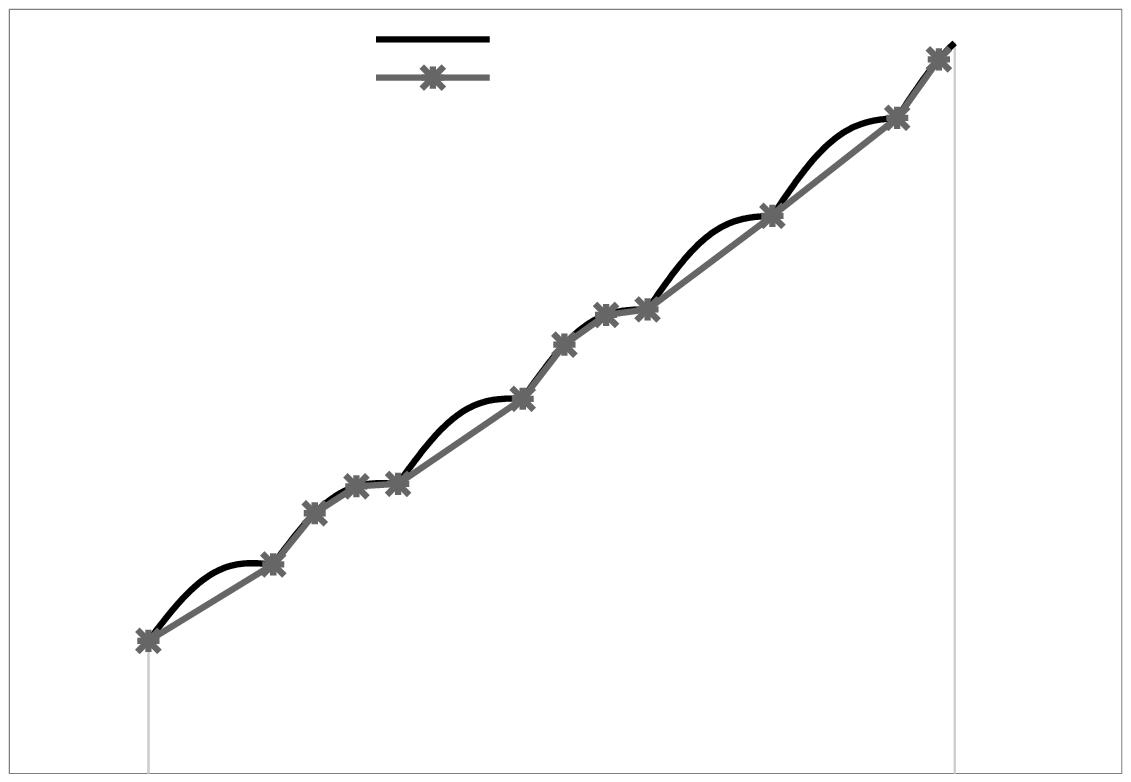}}
\caption{Piecewise-linear fine approximation in important areas, rough on the other areas (left). Fuel cost and its possible final representation with a piecewise-linear function (right).}
\label{fig:cost_lowerB}
\label{fig:cost_lowerC}
\end{figure}

For obtaining increased precision, a large number of breakpoints may be required; this is likely to be the limiting step for large instances.  In order to obviate this problem as much as possible, we may have the cost function represented with greater detail in power regions of potential operation, whereas on other areas the representation may be more rough; in the limit, there may be a single straight line joining two valve points (see Figure~\ref{fig:cost_lowerB}).
Later in this section we propose an algorithm for determining where to expand the number of breakpoints.  The final shape of the piecewise-linear approximation may resemble that of Figure~\ref{fig:cost_lowerC}.

Considering a set of units $\mathcal{U}$ and a set of periods $\mathcal{T}$, for a given unit $u\in\mathcal{U}$ and period $t\in\mathcal{T}$, once the coordinates $(X_{utk},Y_{utk})$ for the set of breakpoints $k=0,\ldots,K$ to be considered is known, the lower approximation of the objective function can easily be included in a mixed-integer programming (MIP) linear model by means of a convex combination of these points:

\begin{alignat}{27}
  & p_{ut} = \sum_{k=0}^{K} X_{utk} z_{utk}, \quad && \forall u \in \mathcal{U}, \forall t \in \mathcal{T}, \label{eq:ccX} \\
  & L_{ut} = \sum_{k=0}^{K} Y_{utk} z_{utk}, \quad && \forall u \in \mathcal{U}, \forall t \in \mathcal{T}, \label{eq:ccY} \\
  & \sum_{k=0}^{K} z_{utk} \leq y_{ut},  \quad && \forall u \in \mathcal{U}, \forall t \in \mathcal{T}. \label{eq:ccz}
\end{alignat}

In a standard convex combination model the last equation is written as $\sum_{k=0}^{K} z_{utk} = 1$.  However, in the case a unit is not committed in a given period, both the production level $p_{ut}$ and the linearized fuel cost $L_{ut}$ must be zero.  For taking this into account, in the right hand size of equation (\ref{eq:ccz}) instead of 1 there must be  variable $y_{ut}$, which indicates whether unit $u$ was committed to produce in period $t$ ($y_{ut} = 1$) or not ($y_{ut}=0$).  Additionally, for the model to be correct, there must be points corresponding to minimum and maximum operating powers, \ie, $X_{ut0} = \Pmin_u$ and $X_{utK} = \Pmax_u$, for each unit $u$ and period $t$.

For the minimization of non-convex functions represented by piecewise-linear segments, as in the present case, it is necessary to introduce binary variables, each corresponding to a $z_{utk}$, limiting which may be non-zero; also an additional constraint for forcing the convex combination to take on two consecutive points is necessary (\ie, there can be $z_{utk}>0$ and $z_{u,t,k+1}>0$ only for one value of $k\in\{0,\ldots,K-1\}$).  Alternatively, we may use a convenience of most modern MIP solvers that assures the same thing by declaring that, for each $u$ and $t$, the set of variables $z_{utk}, \forall k$ forms a so-called  \emph{special ordered set} constraint of type II (SOS2)~\cite{gurobi}.

The objective is to minimize the total production costs,
\begin{alignat}{27}
  \minim \sum_{t\in\mathcal{T}}\sum_{u\in\mathcal{U}} \left(L_{ut} + S_{ut}\right),  \label{eq:objective} 
\end{alignat}
which include start-up costs $S_{ut}$; these are modeled as
\begin{alignat}{27}
  S_{ut} = a_u^{\hot}  s^{\hot}_{ut} + a_u^{\cold}  s_{ut}^{\cold},  \label{eq:startup}\\
\end{alignat}
where $a_u^{\hot}$ is the hot start up cost and variable $s^{\hot}_{ut}=1$ if there was a hot start for unit $u$ in period $t$, 0 otherwise.  Equivalently for cold start up, with $a_u^{\cold}$ and $s_{ut}^{\cold}$.  (See also constraints (\ref{eq:incurred_start_cost}) to~(\ref{eq:switch_off}).)

\subsection{Constraints}
\label{sec:constraints}

In this problem the following constraints will be considered: system power balance and generation limits (\ie, load dispatch), system reserve requirements, unit initial conditions, and unit minimum up and down times; this MIP model was introduced in~\cite{viana2013}.

Demand satisfaction is modeled by constraint~(\ref{eq:load}), and power reserve requirements by constraint~(\ref{eq:reserve}).
\begin{alignat}{27}
  & \sum_{u \in \mathcal{U}} p_{ut} = D_t, \quad && \forall t \in \mathcal{T}, \label{eq:load}\\
  & \sum_{u \in \mathcal{U}} \Pmax_{u} y_{ut} \geq D_t + R_t, \quad && \forall t \in \mathcal{T}. \label{eq:reserve}
\end{alignat}

Power production levels of thermal power units are within the range defined by the technical minimum and maximum production levels in~(\ref{eq:prod_limits}). 
\begin{alignat}{27}
  & \Pmin_u y_{ut} \leq p_{ut} \leq \Pmax_{u} y_{ut}, \quad && \forall u \in \mathcal{U}, \forall t \in \mathcal{T}. \label{eq:prod_limits}
\end{alignat}

When a unit $u$ is switched on, it must remain on for at least $T_u^{\son}$ consecutive periods; similarly, it must be kept off for at least $T_u^{\soff}$ after being switched off.  Constraints (\ref{eq:min_up_ini}) and (\ref{eq:min_down_ini}) model this aspect for the initial state, while constraints (\ref{eq:min_up}) and (\ref{eq:min_down}) do the same for the remaining planning horizon.  In (\ref{eq:min_up_ini}) $\theta^{\son}_u$ represents $\max(0,T_u^{\son} - t^{\prev}_u)$, and $\theta^{\soff}_u$ in (\ref{eq:min_down_ini}) stands for $\max(0,T_u^{\soff} - t^{\prev}_u)$; the previous state of unit $u$ is the parameter $y^{\prev}_u$, which is 1 if the unit was on, 0 if it was off.
\begin{alignat}{27}
  & y_{ut} = 1, \quad && \forall u\in\mathcal{U}: y^{\prev}_u = 1, \for t = 0,\ldots,\theta^{\son}_u,  \label{eq:min_up_ini} \\
  & y_{ut} = 0, \quad && \forall u\in\mathcal{U}: y^{\prev}_u = 0, \for t = 0,\ldots,\theta^{\soff}_u.  \label{eq:min_down_ini}
\end{alignat}
Constraints (\ref{eq:min_up}) and (\ref{eq:min_down}) determine if a unit $u$ is started/switched off in period $t$ ($\xon_{ut}=1, \xoff_{ut}=1$, respectively, or $0$ otherwise); variables $\tau^{\son}_{ut}$ and $\tau^{\soff}_{ut}$ stand for $\max(t - T_u^{\son} + 1, 1)$ and $\max(t - T_u^{\soff} + 1, 1)$, respectively.  
\begin{alignat}{27}
  & \sum_{i = \tau^{\son}_{ut}}^{t} \xon_{ui} \leq  y_{ut},       \quad && \forall u \in \mathcal{U}, \forall t \in \mathcal{T}, \label{eq:min_up} \\
  & \sum_{i = \tau^{\soff}_{ut}}^{t} \xoff_{ui} \leq  1 - y_{ut}, \quad && \forall u \in \mathcal{U}, \forall t \in \mathcal{T}. \label{eq:min_down}
\end{alignat}


Constraints (\ref{eq:incurred_start_cost}) state that every time a unit is switched on, a start-up cost will be incurred.  
\begin{alignat}{27}
  & s^{\hot}_{ut} + s^{\cold}_{ut} = \xon_{ut},  \quad && \forall u \in \mathcal{U}, \forall t \in \mathcal{T}. \label{eq:incurred_start_cost}
\end{alignat}

Constraints (\ref{eq:set_start_cost}) determine the start-up type of each unit, \ie, decide whether it is a cold or a hot start type. It will be a cold start if the unit remained off for more than $ t^{\cold}_u$ periods of time, and a hot start otherwise.

\begin{alignat}{27}
  & y_{ut} - \sum_{i = t - t_u^{\cold} - 1}^{t-1} y_{ui} \leq s^{\cold}_{ut}, \quad && \forall u \in \mathcal{U}, \forall t \in \mathcal{T}. \label{eq:set_start_cost}
\end{alignat}

Constraints (\ref{switch_on}) determine each unit's switch-on variables, and (\ref{eq:switch_off}) determine the switch-off variables.
\begin{alignat}{27}
  & y_{ut} - y_{u,t-1} \leq \xon_{ut},           \quad && \forall u \in \mathcal{U}, \forall t \in \mathcal{T}, \label{switch_on} \\
  & \xoff_{ut} = \xon_{ut} + y_{u,t-1} - y_{ut}, \quad && \forall u \in \mathcal{U}, \forall t \in \mathcal{T}. \label{eq:switch_off}
\end{alignat}


\subsection{Solution approach}
\label{sec:solution}

There is not, to the best of our knowledge, exact solution method for the load dispatching problem when the valve-point loading effect is taken into account, let alone the full unit commitment problem.  The solution approach that we propose achieves this, in the sense that it converges to a solution with an arbitrary degree of precision.

The method is based in the replacement of the function defined by Equation~(\ref{eq:valve}) by a piecewise-linear function; a careful selection of the breakpoints allows the construction of a model that provides a lower bound to the exact function.  Upon a feasible solution obtained by this linear model, we can evaluate the true function; this will provide an upper bound to the value of the objective.  As described in Algorithm~\ref{alg:solution}, if the deviation between the upper and lower bounds are within a user-specified tolerance, then the method will stop.  Otherwise, the number of breakpoints is increased (though only in the required intervals), thus leading to a better approximation, and the optimization problem is resolved.  This process is repeated until the relative deviation between upper and lower bounds is small enough (or the allowed CPU time is exceeded).  Notice that the solution of an iteration is feasible to the problem of the next, the only changes are in its evaluation.  We can, therefore, supply it as a starting point to the MIP solver in line (\ref{itm:solve}) of the algorithm.

\begin{myalg}[htb!]
\caption{Solution procedure.}
\label{alg:solution}
\begin{footnotesize}
\alginout{UCP instance, CPU limit, tolerance}{A solution to the input instance}
\begin{algtab}
  \algforeach{unit $u \in \mathcal{U}$, period $t \in \mathcal{T}$:}
  $B_{ut} := $ set of breakpoints containing only the valve points\\
  \algend
  \algwhile{\algtrue}
  solve problem with piecewise-linear objective, breakpoints $B$  \label{itm:solve}\\
  \algif{$UB - LB <$ tolerance \algor CPU limit exceeded:}
  \algreturn solution\\
  \algend
  \algforeach{$u, t$ with coarse objective representation in the solution:}
  divide the inter-valve interval containing $p_{ut}$ into $K$ segments\\
  add corresponding points to $B_{ut}$\\
  \algend
  \algif{there were no $u,t$ with coarse objective in $p_{ut}$:}
  $K := 2 K$\\
  for each $u,t$ recompute $B_{ut}$ based on $K$ segments in the inter-valve interval containing $p_{ut}$\\
  \algend
  \algend
\end{algtab}  
\end{footnotesize}
\end{myalg}

The value $K$ used in this algorithm should ideally be setup in such a way that the initial approximation is just good enough to solve the instance with sufficient precision.  If $K$ is set too high, there will be too many intervals in the piecewise-linear representation of the objective function, leading to great precision but requiring unnecessary CPU for solving the problem; if is it set too low, additional iterations will be necessary for providing the required precision (each iteration doubling the value of $K$), and again unnecessary CPU will be used.

The most important properties of this algorithm are summarized next.  Let us denominate $P_1$ the problem of maximizing the true objective,
\begin{alignat}{27}
  \minim \sum_{t\in\mathcal{T}}\sum_{u\in\mathcal{U}}
  \left[ S_{ut} + y_{ut} F_{ut}(p_{ut}) \right] \label{eq:trueobj}
\end{alignat}
where $$F_{ut}(p) = a_u + b_u p + c_u p^2 + \left|e_u \sin\left(f_u (\Pmin_{ut}-p) \right) \right|,$$
under constraints (\ref{eq:startup}) to~(\ref{eq:switch_off}), and $P_2$ its linearized counterpart, \ie, the problem of maximizing~(\ref{eq:objective}) under the same constraints and additionally (\ref{eq:ccX}), (\ref{eq:ccY}) and~(\ref{eq:ccz}).

\begin{figure}[!htbp]
\centering
\begin{minipage}{1.07\linewidth}
\hspace{-.05\linewidth}
\resizebox{.51\linewidth}{!}{\input{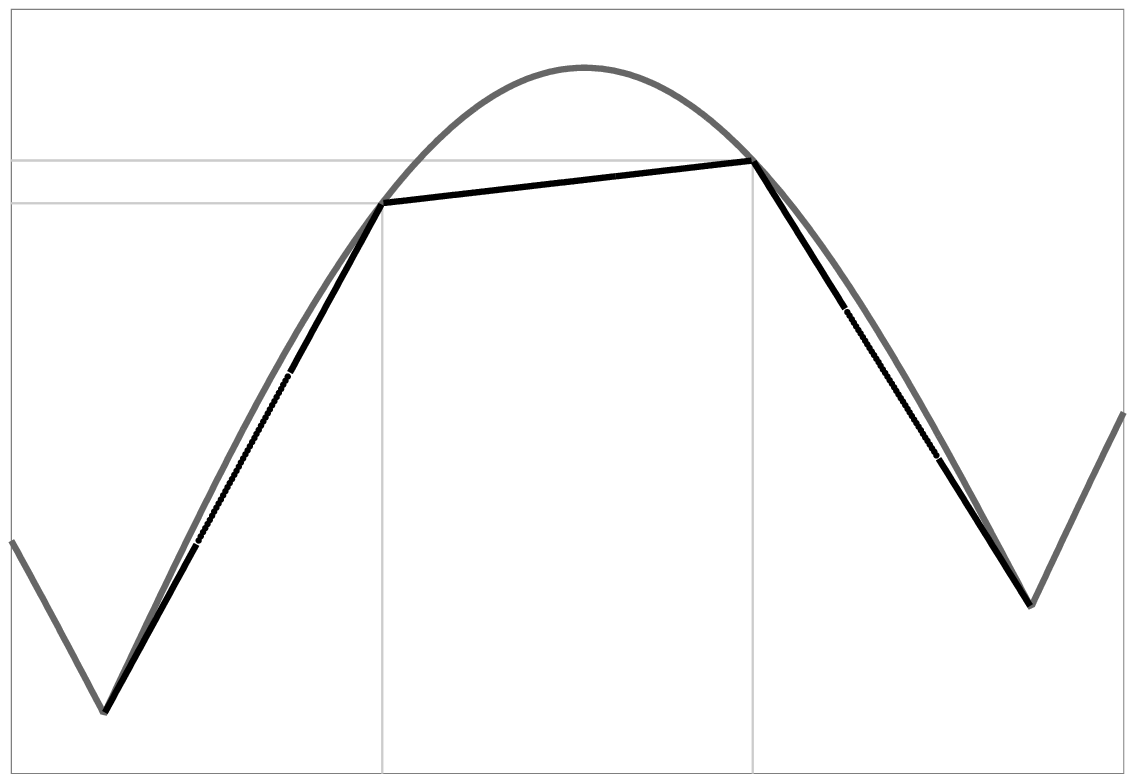}}
\resizebox{.51\linewidth}{!}{\input{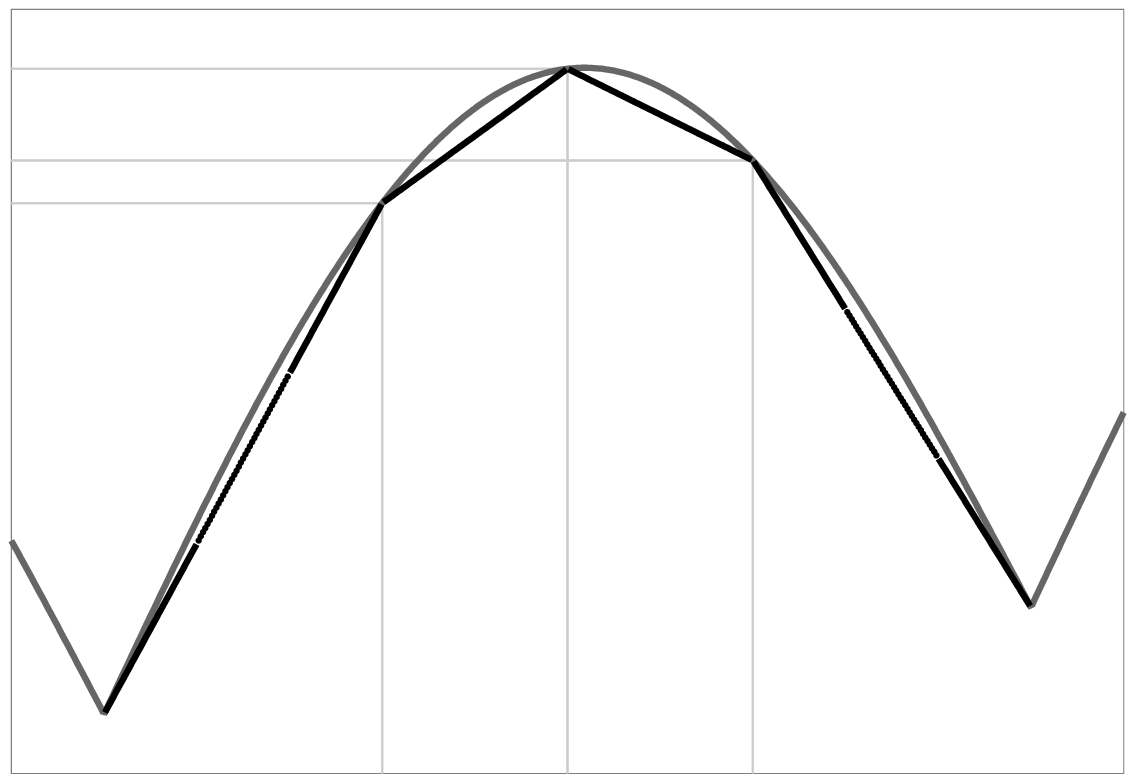}}
\end{minipage}
\caption{Linear approximations in consecutive iterations.}
\label{fig:cost_shape}
\end{figure}

\begin{property}
  \label{prop:feas}
  Feasible solutions for $P_2$ are feasible for $P_1$.
\end{property}
\begin{proof}
  $P_2$ includes all the variables of $P_1$, and contraints of $P_2$ are a superset of constraints of $P_1$.
\end{proof}

\begin{property}
  The evaluation of the objective function~(\ref{eq:trueobj}) at a feasible solution obtained for $P_2$ is an upper bound to the optimum of $P_1$.
\end{property}
\begin{proof}
  Follows from Property~\ref{prop:feas} and from the notion of optimum: the objective at any feasible solution of a minimization problem provides an upper bound to its optimum.
\end{proof}

\begin{property}
  An optimum of $P_2$ is a lower bound to the optimum of $P_1$.
\end{property}
\begin{proof}
  Between two consecutive valve points the function $F(p)$ defined by Equation~\ref{eq:valve} is concave as long as $c \leq e f^2 / 2$ (see~\ref{sec:concavity}); therefore, if this condition is observed, any piecewise-linear interpolation between consecutive valve points is not larger than $F$.  In the construction of the piecewise-linear approximation we are considering all the valve points, as well as the extreme point of the domain of $F$, as interpolation points; hence, the linear interpolation is not larger than $F$ in its domain.  We can thus conclude that the minimum of $P_2$ --- where the objective is a linear interpolation in the previous conditions --- is not larger than the minimum of $P_1$, where the objective is function~$F$.
\end{proof}

\begin{property}
  With tolerance $\epsilon = 0$, when the number of iterations increases the solution provided by Algorithm~\ref{alg:solution} converges to the optimum solution of~$P_1$.
\end{property}
\begin{proof}
  Let us denote the linear approximation of $F$ at an iteration $i$ by $G_i$, and define $G = \lim_{i\rightarrow\infty}G_i$.  Note that $G_j(p) \geq G_i(p)$, for any $j>i$, and for all $p$ in the domain of $F$.  As functions $G_i$ are non-decreasing with increasing $i$ and are limited above by $F$, the sequence $G_i$ converges.  

It remains to prove that the optimum solution with $G_i$ as the objective function, when $i$ tends to infinity, is an optimum solution with $F$ as the objective function.  We will assume that this is not true, \ie, that the optimum $F^*$ is strictly greater than $G^*$, and show that this leads to a contradiction.  Notice that $F^* > G^*$ is equivalent to $F^* - G^* = \epsilon^*$, with $\epsilon^*>0$, which implies that there must be some $u,t$ such that $F_{ut}(p^*) - G_{ut}(p^*) = \epsilon' > 0$.  However, lines (11) and (12) of Algorithm~\ref{alg:solution} imply that the number of segments in each inter-valve interval containing $p_{ut}^*$ tends to infinity when the number of iterations tends to infinity.  This means that $G_{ut}(p^*)$ becomes arbitrarily close to $F_{ut}(p^*)$, and hence $\epsilon'$ cannot be positive.
\end{proof}

\section{Computational results}
\label{sec:results}

In order to run economic load dispatch instances with the current model we must set up an instance with only one period; generators must have the previous status as operating, and must be forced to keep operating by adjusting the minimum number of on periods.  We have thus prepared data for the most widely used benchmark instances for load dispatching with valve-point loading effect: \emph{eld13} (first described in \cite{wood1984}) and \emph{eld40} (first described in \cite{sinha2003}); we have used the tables available in \cite{coelho2006}\footnote{All the data and programs used are available in this paper's Internet page, \url{http://www.dcc.fc.up.pt/~jpp/code/valve}}.

The setup used was the following: a computer with an Intel Xeon processor at 3.0 GHz and running Linux version 2.6.32, using Gurobi version 5.0.1~\cite{gurobi}.  Only one thread was assigned to this experiment.  Models were written in the Python language and the default Gurobi parameters were used until reaching the user-defined tolerance; then, parameters were changed to the maximum precision supported by Gurobi.  The aim of this procedure is to avoid spending too much CPU in the initial iterations, where the representation of the objective function is still rather course.

\begin{table}[!htbp]
  \centering
  \begin{tabular}{l|r|r|r|r|r|r}
    \hline
    Instance & Load      & Lower bound & Upper bound & Error    & CPU time      & $N$   \\\hline
    \emph{eld13}    & 1800      &  17963.83   &  17963.83   & $<$1.e-7 & 19.1s  & 34  \\
    \emph{eld13}    & 2520      &  24169.92   &  24169.92   & $<$1.e-7 &  1.8s  & 19  \\
    \emph{eld40}    & 10500     & 121412.53   & 121412.54   & $<$1.e-7 &  7.6s  & 27  \\
    \hline
  \end{tabular}
  \caption{Solution to economic load dispatch instances: bounds to the objective value, maximum relative error, CPU time employed, and number of iterations performed.}
  \label{tab:eld}
\end{table}

The results presented in Table~\ref{tab:eld} are optima, and are approximately equal to the best solutions found by the best methods available in the literature~\cite{SrinivasaReddy2013342}.\footnote{Results presented in~\cite{Coelho20102580} are better than the optimum for \emph{eld13}, indicating that there is likely an error in that paper.  Possibly this is due to a rounding error; if we calculate the cost for the solutions provided in the paper we obtain values greater than those reported.}

\begin{figure}[!t]
\begin{minipage}{1.2\linewidth}
\begin{centering}
~\hspace{-0.2\textwidth}
\includegraphics[width=.48\columnwidth]{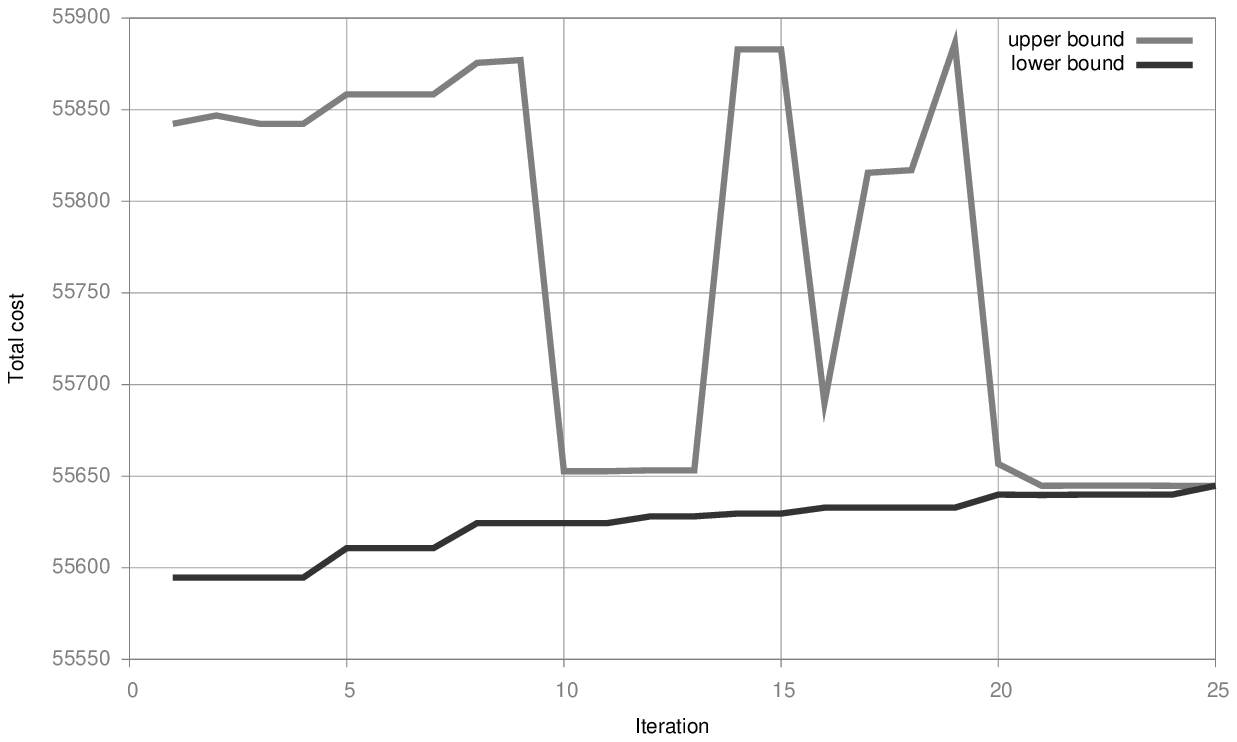}
\includegraphics[width=.48\columnwidth]{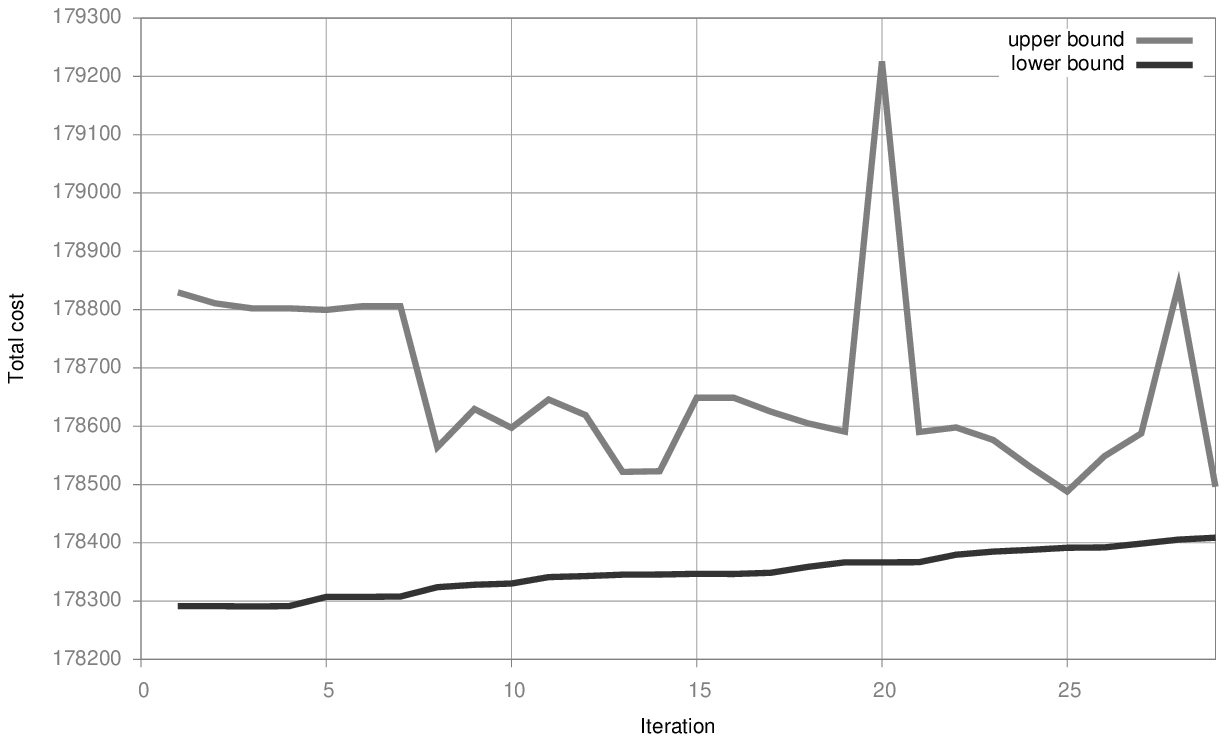}
\par\end{centering}
\end{minipage}
\caption{\textit{Bounds as a function of the iteration number, for instance \emph{ucp40} with one period (optimal solution was found) and three periods (final solution was not optimal).}\label{fig:bound}}
\end{figure}

As for the unit commitment problem, to the best of our knowledge there are no test instances taking into account the valve-point effect.  We have therefore prepared a set of benchmarks, based on the instances for unit commitment provided in~\cite{kazarlis1996}, and on the instances for economic load dispatch mentioned above.  Instance \emph{ucp10} corresponds to the standard instance described in~\cite{kazarlis1996}; \emph{ucp5} is based on that instance, by selecting only half of the units and adjusting demand.  Instances \emph{ucp13} and \emph{ucp40} have costs taken from \emph{eld13} and \emph{eld40}, respectively, and data concerning multi-period operation based on~\cite{kazarlis1996}.  Notice that the difference between instances \emph{eld13} and \emph{ucp13} with one period (besides possible differences in the demand value) is that in the former the commitment decision has already been made, and all units considered must be operating; the same for \emph{eld40} and \emph{ucp40}.  These results are presented in Table~\ref{tab:ucp}.  Figure~\ref{fig:bound} plots the evolution, for a particular case, of the lower and upper bounds with respect to the iteration number; as expected, the lower bound is non-decreasing, whereas the upper bound may increase (when the linear approximation at the solution of the previous iteration was poor).

\begin{table}[!htbp]
  \centering
  \begin{tabular}{l|r|r|r|r|r@{.}l|r}
    \hline
    Instance & Periods   & Lower bound & Upper bound & Error  & \multicolumn{2}{r|}{CPU time} & $N$ \\\hline
\emph{ucp5}  & 1  &     7352.46 &    7352.46 & $<$1.e-7   &    0&12  &  5   \\ 
\emph{ucp5}  & 3  &    24506.12 &   24506.12 & $<$1.e-7   &    0&29  &  5   \\ 
\emph{ucp5}  & 6  &    59127.38 &   59127.39 & $<$1.e-7   &    1&2   &  6   \\ 
\emph{ucp5}  & 12 &   153575.68 &  153575.69 & $<$1.e-7   &   54&    &  7   \\ 
\emph{ucp5}  & 24 &   308482.03 &  309152.47 & 0.22\%     & 3600&    &  3   \\ \hline 
\emph{ucp10} & 1  &    13826.85 &   13826.85 & $<$1.e-7   &    0&17  &  6   \\ 
\emph{ucp10} & 3  &    45929.36 &   45929.36 & $<$1.e-7   &    1&3   &  7   \\ 
\emph{ucp10} & 6  &   109919.97 &  109919.98 & $<$1.e-7   &   31&    &  9   \\ 
\emph{ucp10} & 12 &   285808.54 &  286067.44 & 0.091\%    & 3600&    &  5   \\ 
\emph{ucp10} & 24 &   572943.80 &  574554.02 & 0.28\%     & 3600&    &  3   \\ \hline 
\emph{ucp13} & 1  &    11701.28 &   11701.28 & $<$1.e-7   &    0&15  &  6   \\ 
\emph{ucp13} & 3  &    38849.84 &   38849.84 & $<$1.e-7   &    3&3   & 11   \\ 
\emph{ucp13} & 6  &    91405.97 &   91783.98 & 0.41\%     & 3600&    & 10   \\ 
\emph{ucp13} & 12 &   231587.47 &  232537.37 & 0.41\%     & 3600&    &  4   \\ 
\emph{ucp13} & 24 &   464053.01 &  466186.98 & 0.46\%     & 3600&    &  3   \\ \hline 
\emph{ucp40} & 1  &    55644.79 &   55644.79 & $<$1.e-7   &    3&9   & 25   \\ 
\emph{ucp40} & 3  &   178395.51 &  178547.12 & 0.085\%    & 3600&    & 30   \\ 
\emph{ucp40} & 6  &   416108.42 &  416605.93 & 0.12\%     & 3600&    & 15   \\ 
\emph{ucp40} & 12 &  1112370.94 & 1113801.09 & 0.13\%     & 3600&    & 11   \\ 
\emph{ucp40} & 24 &  2235970.56 & 2238504.45 & 0.11\%     & 3600&    &  7   \\ 
    \hline
  \end{tabular}
  \caption{Solution to unit commitment instances: bounds to the objective value, maximum relative error, CPU time employed, and number of iterations performed.}
  \label{tab:ucp}
\end{table}

\section{Conclusions}
\label{sec:conclusions}

The main contribution of this paper is a model and an algorithm for solving the unit commitment problem taking into consideration the valve-point loading effect in fuel cost computation.  This model can be used for planning which units to commit in each period so as to match demand, meeting reserve constraints, and under a set of technological constraints.  The algorithm iteratively calls a general-purpose mixed-integer programming solver for tackling optimization subproblems, and uses the solution to improve the accuracy of fuel cost approximation until reaching an error below a given tolerance, or reaching a given limit on CPU usage.

The model can also be used for load dispatching, when only one period is considered and the subset of units that will operate is already selected; this is known as the economic load dispatch problem with valve-point effect; the method proposed was able to solve well known benchmark instances of this problem for which the optimum was not known.

To the best of our knowledge, there are no standard benchmarks for unit commitment taking into account the valve-point loading effect.  We present a set of them, based on data from load dispatch and unit commitment.  Using our algorithm for their solution, the optimum was obtained for a subset of instances, and very good approximations, with a relative error below 1\%, were found for the remaining.

\section*{Acknowledgements}
This work was partly funded by project ``NORTE-07-0124-FEDER-000057'', under the North Portugal Regional Operational Programme (ON.2 – O Novo Norte) and the National Strategic Reference Framework through the European Regional Development Fund, and by national funds through Funda{\c c}{\~a}o para a Ci{\^e}ncia e a Tecnologia (FCT), and by FCT (under Project PTDC/EGE-GES/099120/2008) through ``Programa Operacional Tem\'atico Factores de Competitividade (COMPETE)'' of ``Quadro Comunit\'ario de Apoio III'', partially funded by FEDER.

\bibliographystyle{model3-num-names}
\bibliography{ucp-valve}

\begin{thebibliography}{17}
\providecommand{\natexlab}[1]{#1}
\providecommand{\url}[1]{\texttt{#1}}
\providecommand{\urlprefix}{URL }
\expandafter\ifx\csname urlstyle\endcsname\relax
  \providecommand{\doi}[1]{doi:\discretionary{}{}{}#1}\else
  \providecommand{\doi}{doi:\discretionary{}{}{}\begingroup
  \urlstyle{rm}\Url}\fi
\providecommand{\eprint}[2][]{\url{#2}}
\providecommand{\BIBand}{and}
\providecommand{\bibinfo}[2]{#2}
\ifx\xfnm\undefined \def\xfnm[#1]{\unskip,\space#1}\fi
\bibitem[{Al-Sumait et~al.(2007)Al-Sumait, Al-Othman and
  Sykulski}]{AlSumait2007720}
\bibinfo{author}{Al-Sumait\xfnm[ J.]}, \bibinfo{author}{Al-Othman\xfnm[ A.]},
  \bibinfo{author}{Sykulski\xfnm[ J.]}.
\newblock \bibinfo{title}{Application of pattern search method to power system
  valve-point economic load dispatch}.
\newblock \bibinfo{journal}{International Journal of Electrical Power \& Energy
  Systems}
  \bibinfo{year}{2007};\bibinfo{volume}{29}(\bibinfo{number}{10}):\bibinfo{pages}{720
  -- 730}.
\bibitem[{Hemamalini and Simon(2011)}]{Hemamalini2011868}
\bibinfo{author}{Hemamalini\xfnm[ S.]}, \bibinfo{author}{Simon\xfnm[ S.P.]}.
\newblock \bibinfo{title}{Dynamic economic dispatch using artificial immune
  system for units with valve-point effect}.
\newblock \bibinfo{journal}{International Journal of Electrical Power \& Energy
  Systems}
  \bibinfo{year}{2011};\bibinfo{volume}{33}(\bibinfo{number}{4}):\bibinfo{pages}{868
  -- 874}.
\bibitem[{Roy et~al.(2013)Roy, Roy and Chakrabarti}]{Roy20134244}
\bibinfo{author}{Roy\xfnm[ P.]}, \bibinfo{author}{Roy\xfnm[ P.]},
  \bibinfo{author}{Chakrabarti\xfnm[ A.]}.
\newblock \bibinfo{title}{Modified shuffled frog leaping algorithm with genetic
  algorithm crossover for solving economic load dispatch problem with
  valve-point effect}.
\newblock \bibinfo{journal}{Applied Soft Computing}
  \bibinfo{year}{2013};\bibinfo{volume}{13}(\bibinfo{number}{11}):\bibinfo{pages}{4244
  -- 4252}.
\bibitem[{Reddy and Vaisakh(2013)}]{SrinivasaReddy2013342}
\bibinfo{author}{Reddy\xfnm[ A.S.]}, \bibinfo{author}{Vaisakh\xfnm[ K.]}.
\newblock \bibinfo{title}{Shuffled differential evolution for economic dispatch
  with valve point loading effects}.
\newblock \bibinfo{journal}{International Journal of Electrical Power \& Energy
  Systems}
  \bibinfo{year}{2013};\bibinfo{volume}{46}(\bibinfo{number}{0}):\bibinfo{pages}{342
  -- 352}.
\bibitem[{Viana and Pedroso(2013)}]{viana2013}
\bibinfo{author}{Viana\xfnm[ A.]}, \bibinfo{author}{Pedroso\xfnm[ J.P.]}.
\newblock \bibinfo{title}{A new {MILP}-based approach for unit commitment in
  power production planning}.
\newblock \bibinfo{journal}{International Journal of Electrical Power and
  Energy Systems}
  \bibinfo{year}{2013};\bibinfo{volume}{44}:\bibinfo{pages}{997--1005}.
\bibitem[{Sen and Kothari(1998)}]{sen1998}
\bibinfo{author}{Sen\xfnm[ S.]}, \bibinfo{author}{Kothari\xfnm[ D.]}.
\newblock \bibinfo{title}{Optimal thermal generating unit commitment: a
  review}.
\newblock \bibinfo{journal}{International Journal of Electrical Power \& Energy
  Systems}
  \bibinfo{year}{1998};\bibinfo{volume}{20}(\bibinfo{number}{7}):\bibinfo{pages}{443--451}.
\bibitem[{Padhy(2004)}]{PAD04}
\bibinfo{author}{Padhy\xfnm[ N.P.]}.
\newblock \bibinfo{title}{Unit commitment -- a bibliographical survey}.
\newblock \bibinfo{journal}{IEEE Transactions in Power Systems}
  \bibinfo{year}{2004};\bibinfo{volume}{19}(\bibinfo{number}{2}):\bibinfo{pages}{1196--1205}.
\bibitem[{Yamin(2004)}]{YAM04}
\bibinfo{author}{Yamin\xfnm[ H.Y.]}.
\newblock \bibinfo{title}{Review on methods of generation scheduling in
  electric power systems}.
\newblock \bibinfo{journal}{Electric Power Systems Research}
  \bibinfo{year}{2004};\bibinfo{volume}{69}:\bibinfo{pages}{227--248}.
\bibitem[{L{\'o}pez et~al.(2012)L{\'o}pez, Ceciliano-Meza, Moya and
  G{\'o}mez}]{Alvarez2012}
\bibinfo{author}{L{\'o}pez\xfnm[ J.{\'A}.]},
  \bibinfo{author}{Ceciliano-Meza\xfnm[ J.L.]}, \bibinfo{author}{Moya\xfnm[
  I.G.]}, \bibinfo{author}{G{\'o}mez\xfnm[ R.N.]}.
\newblock \bibinfo{title}{A {MIQCP} formulation to solve the unit commitment
  problem for large-scale power systems}.
\newblock \bibinfo{journal}{International Journal of Electrical Power \& Energy
  Systems}
  \bibinfo{year}{2012};\bibinfo{volume}{36}(\bibinfo{number}{1}):\bibinfo{pages}{68
  -- 75}.
\bibitem[{Simoglou et~al.(2010)Simoglou, Biskas and Bakirtzis}]{Simoglou2010}
\bibinfo{author}{Simoglou\xfnm[ C.K.]}, \bibinfo{author}{Biskas\xfnm[ P.N.]},
  \bibinfo{author}{Bakirtzis\xfnm[ A.G.]}.
\newblock \bibinfo{title}{Optimal self-scheduling of a thermal producer in
  short-term electricity markets by milp}.
\newblock \bibinfo{journal}{IEEE Transactions on Power Systems}
  \bibinfo{year}{2010};\bibinfo{volume}{25}(\bibinfo{number}{4}):\bibinfo{pages}{1965
  -- 1977}.
\bibitem[{Lima et~al.(2013)Lima, Marcovecchio, Novais and Grossmann}]{Lima2013}
\bibinfo{author}{Lima\xfnm[ R.M.]}, \bibinfo{author}{Marcovecchio\xfnm[ M.G.]},
  \bibinfo{author}{Novais\xfnm[ A.Q.]}, \bibinfo{author}{Grossmann\xfnm[
  I.E.]}.
\newblock \bibinfo{title}{On the computational studies of deterministic global
  optimization of head dependent short-term hydro scheduling}.
\newblock \bibinfo{journal}{IEEE Transactions on Power Systems}
  \bibinfo{year}{2013};\bibinfo{volume}{28}(\bibinfo{number}{4}):\bibinfo{pages}{4336
  -- 4347}.
\bibitem[{Kazarlis et~al.(1996)Kazarlis, Bakirtzis and Petridis}]{kazarlis1996}
\bibinfo{author}{Kazarlis\xfnm[ S.A.]}, \bibinfo{author}{Bakirtzis\xfnm[
  A.G.]}, \bibinfo{author}{Petridis\xfnm[ V.]}.
\newblock \bibinfo{title}{A genetic algorithm solution to the unit commitment
  problem}.
\newblock \bibinfo{journal}{IEEE Transactions on Power Systems}
  \bibinfo{year}{1996};\bibinfo{volume}{11}:\bibinfo{pages}{83--92}.
\bibitem[{Wood and Wollenberg(1984)}]{wood1984}
\bibinfo{author}{Wood\xfnm[ A.]}, \bibinfo{author}{Wollenberg\xfnm[ B.]}.
\newblock \bibinfo{title}{Power Generation, Operation, and Control}.
\newblock \bibinfo{address}{New York, NY}: \bibinfo{publisher}{John Wiley \&
  Sons}; \bibinfo{year}{1984}.
\bibitem[{dos Santos~Coelho and Mariani(2006)}]{coelho2006}
\bibinfo{author}{dos Santos~Coelho\xfnm[ L.]}, \bibinfo{author}{Mariani\xfnm[
  V.]}.
\newblock \bibinfo{title}{Combining of chaotic differential evolution and
  quadratic programming for economic dispatch optimization with valve-point
  effect}.
\newblock \bibinfo{journal}{IEEE Transactions on Power Systems}
  \bibinfo{year}{2006};\bibinfo{volume}{21}(\bibinfo{number}{2}):\bibinfo{pages}{989--996}.
\newblock \bibinfo{note}{(See also a correction to this paper)}.
\bibitem[{dos Santos~Coelho and Mariani(2010)}]{Coelho20102580}
\bibinfo{author}{dos Santos~Coelho\xfnm[ L.]}, \bibinfo{author}{Mariani\xfnm[
  V.C.]}.
\newblock \bibinfo{title}{An efficient cultural self-organizing migrating
  strategy for economic dispatch optimization with valve-point effect}.
\newblock \bibinfo{journal}{Energy Conversion and Management}
  \bibinfo{year}{2010};\bibinfo{volume}{51}(\bibinfo{number}{12}):\bibinfo{pages}{2580
  -- 2587}.
\bibitem[{{Gurobi Optimization, Inc.}(2012)}]{gurobi}
\bibinfo{author}{{Gurobi Optimization, Inc.}\xfnm[]}.
\newblock \bibinfo{title}{Gurobi Optimizer Reference Manual, Version 5.0}.
\newblock \bibinfo{address}{http://www.gurobi.com}; \bibinfo{year}{2012}.
\bibitem[{Sinha et~al.(2003)Sinha, Chakrabarti and Chattopadhyay}]{sinha2003}
\bibinfo{author}{Sinha\xfnm[ N.]}, \bibinfo{author}{Chakrabarti\xfnm[ R.]},
  \bibinfo{author}{Chattopadhyay\xfnm[ P.K.]}.
\newblock \bibinfo{title}{Evolutionary programming techniques for economic load
  dispatch}.
\newblock \bibinfo{journal}{IEEE Transactions on Evolutionary Computation}
  \bibinfo{year}{2003};\bibinfo{volume}{7}(\bibinfo{number}{1}):\bibinfo{pages}{83--94}.

\end{thebibliography}

\appendix
\section{Notation}

\subsection*{Constants}

\begin{itemize}
\item $T$ -- length of the planning horizon.
\item $\mathcal{T} = \{1,\ldots,T\}$ -- set of planning periods.
\item $\mathcal{U}$ -- set of units.
\item $\Pmin_u, \Pmax_u$ -- minimum and maximum production levels for unit $u$.
\item $T_u^{\son}, T_u^{\soff}$ -- minimum number of periods that unit $u$ must be kept switched on/off.
\item $D_t$ -- system load requirements in period~$t$.
\item $R_t$ -- spinning reserve requirements in period~$t$.
\item $a_u, b_u, c_u, e_u, f_u$ -- fuel cost parameters for unit~$u$.
\item $a^{\hot}_{u}, a^{\cold}_{u}$ -- hot and cold start up costs for unit~$u$.
\item $t^{\cold}_u$ -- number of periods after which the start up of unit $u$ is evaluated as cold.
\item $y^{\prev}_u$ -- previous state of unit $u$ (1 if on, 0 if off).
\item $t^{\prev}_u$ -- number of periods unit $u$ has been on or off prior to the first period of the planning horizon.
\end{itemize}

\subsection*{Decision variables}
\begin{itemize}
\item $y_{ut}$ -- 1 if unit $u$ is on in period $t$, $0$ otherwise.
\item $p_{ut}$ -- production level of unit $u$, in period $t$.
\end{itemize}

\subsection*{Auxiliary variables}
\begin{itemize}
\item $\xon_{ut}, \xoff_{ut}$ -- 1 if unit $u$ is started/switched off in period $t$, $0$ otherwise;
\item $s^{\hot}_{ut}$ -- 1 if unit $u$ has a hot start in period $t$, $0$ otherwise;
\item $s^{\cold}_{ut}$ -- 1 if unit $u$ has a cold start in period $t$, $0$ otherwise;
\end{itemize}

\subsection*{Production costs}
\begin{itemize}
\item $F_{ut}$ -- fuel cost for unit $u$ in period  $t$.
\item $S_{ut}$ -- start up cost for unit $u$ in period $t$.
\end{itemize}

\section{Conditions for concavity of the fuel costs}
\label{sec:concavity}

In this appendix we state the conditions for concavity of the fuel cost function $F(p) = a + b p + c p^2 + \left|e \sin\left(f (\Pmin-p)\right)\right|$ (Equation~\ref{eq:valve}) between two valve points.
Recall that a differentiable function is concave in a given region if its second derivative is less than or equal to zero in all points of that region.


If the term within the absolute value, $e \sin\left(f (\Pmin - p)\right)$, is positive, the first derivative of $F(p)$ is $F'(p) = 2 c p + b - e f \cos\left(f (\Pmin - p)\right)$, and the second derivative is 
$F''(p) = 2c - e f^2 \sin(f (\Pmin - p)).$
For the function $F$ to be concave there must be $F''(p) \leq 0$, \ie, $c \leq \frac{e f^2}{2} \sin\left(f (\Pmin - p)\right)$.  As  $\sin(x) \leq 1, \forall x$, $F$ is concave for $c \leq e f^2 / 2$.

Similarly, if the term within the absolute value is negative, the first derivative of $F(p)$ is $F'(p) = 2 c p + b + e f \cos\left(f (\Pmin - p)\right)$, and the second derivative is 
$F''(p) = 2 c + e f^2 \sin \left(f (\Pmin - p)\right).$
Again, for $F$ to be concave there must be $c \leq$ $- \frac{e f^2}{2} \sin\left(f (P^{min} - p)\right)$.  As $\sin(x) \geq -1, \forall x$, this condition is observed for $c \leq e f^2 / 2$.

Therefore, for  
$$c \leq e f^2 / 2$$
the function $F(p)$ is concave between two valve points. This condition is valid for all the instances analyzed in this paper.

\end{document}